\documentclass{amsart}
\usepackage{amsmath,amssymb,amsfonts,enumerate,amsthm}
\usepackage[pdftitle={Title of Article}, verbose,
            colorlinks=true,
            naturalnames=true,
            linkcolor=blue,
            citecolor=red,
            pdfstartview=XYZ]{hyperref}

\newcommand{\Hom}{\mbox{Hom}\,}
\newcommand{\Ext}{\mbox{Ext}\,}
\newcommand{\Tor}{\mbox{Tor}\,}
\newcommand{\Spec}{\mbox{Spec}\,}

\newcommand{\Supp}{\mbox{Supp}\,}

\newcommand{\depth}{\mbox{depth}\,}
\newcommand{\width}{\mbox{width}\,}
\renewcommand{\dim}{\mbox{dim}\,}

\newcommand{\pd}{\mbox{proj.dim}\,}
\newcommand{\id}{\mbox{inj.dim}\,}

\newcommand{\gid}{\mbox{Gid}\,}
\newcommand{\Gid}{\mbox{Gid}\,}

\renewcommand{\H}{\mbox{H}}

\newcommand{\lo}{\longrightarrow}

\newcommand{\fm}{\mathfrak{m}}
\newcommand{\m}{\mathfrak{m}}
\newcommand{\fp}{\mathfrak{p}}

\newcommand{\fq}{\mathfrak{q}}
\newcommand{\q}{\mathfrak{q}}
\newcommand{\fn}{\mathfrak{n}}
\newcommand{\n}{\mathfrak{n}}
\newcommand{\varph}{\varphi}
\newcommand{\otime}{\otimes}
\newtheorem{thm}{Theorem}[section]
\newtheorem{cor}[thm]{Corollary}
\newtheorem{lem}[thm]{Lemma}
\newtheorem{prop}[thm]{Proposition}
\newtheorem{defn}[thm]{Definition}

\numberwithin{equation}{section}

\begin{document}
\bibliographystyle{amsplain}

\title[Local Cohomology \& Gorenstein inj. dim. over local
homomorphisms]{Local Cohomology and Gorenstein injective dimension over local
homomorphisms}

\author{Leila Khatami}
\address{Leila Khatami\\ The Abdus Salam ICTP, Strada Costiera 11, 34100
Trieste, Italy}

\author{Massoud Tousi}
\address{Massoud Tousi\\Department of Mathematics, Shahid Beheshti
University, Tehran, Iran,\\ and Institute for Theoretical Physics
and Mathematics (IPM).}

\author{Siamak Yassemi}


\address{Department of Mathematics,
University of Tehran, Tehran, Iran,\\ and Institute for Theoretical
Physics and Mathematics (IPM).}

\thanks{Emails:  \href{mailto:lkhatami@ictp.it}{lkhatami@ictp.it};
\href{mailto:mtousi@ipm.ir}{mtousi@ipm.ir};
\href{mailto:yassemi@ipm.ir}{yassemi@ipm.ir}}

\keywords{Local cohomology, Gorenstein injective dimension, Bass'
formula}

\subjclass[2000]{13D05, 13D45, 14B15}

\begin{abstract}

\noindent Let $\varphi: (R,\fm) \to (S,\fn)$ be a local
homomorphism of commutative noetherian local rings. Suppose that
$M$ is a finitely generated $S$-module. A generalization of
Grothendieck's non-vanishing theorem is proved for $M$ (i.e. the
Krull dimension of $M$ over $R$ is the greatest integer $i$ for
which the $i$th local cohomology module of $M$ with respect to
$\fm$, $H^i_\fm(M)$, is non-zero). It is also proved that the
Gorenstein injective dimension of $M$, if finite, is bounded below
by dimension of $M$ over $R$ and is equal to the supremum of
$\depth R_\fp$, where $\fp$ runs over the support of $M$ as an
$R$-module.
\end{abstract}

\maketitle

\section*{Introduction}

Let $R$ be a non-trivial commutative noetherian local ring. An
$R$-module $M$ is called {\it finite over a local homomorphism} if
there exists a local homomorphism of noetherian local rings $R\to
S$ such that $M$ is a finite (that is finitely-generated)
$S$-module and the $S$-action is compatible with the action of
$R$. Studied by Apassov, Avramov, Christensen, Foxby, Iyengar,
Miller, Sather-Wagstaff and others, cf. \cite{Ap,AF,AIM,CI,IS},
homological properties of finite modules over (local)
homomorphisms are shown to extend the those of finite modules.

It is well-known that over a local ring, injective dimension of a
finite module is either infinite or equals the depth of the ring
and is not less than the dimension of the module. This is known as
the {\it Bass formula}.

A similar formula is also showed to be true for a module finite
over a local homomorphism (cf. \cite{AIM,TY,KY}).

In this paper, we deal with {\it Gorenstein injective} dimension
of a module which is finite over a local homomorphism. Introduced
by Enochs and Jenda \cite{EJ1,EJ2}, Gorenstein
injective dimension is the dual notion to the G-dimension, due to
Auslander and Bridger \cite{ABr}.

In \cite{KY} and \cite{Y}, a generalized Bass formula is proved
for finite modules of finite Gorenstein injective dimension.\\
{\bf Theorem.} Let $R$ be a noetherian local ring and $M$ be a
finite module of finite Gorenstein injective dimension. Then
$$\dim_R M \leq \Gid_R(M)=\sup \{\depth R_\fp \, | \, \fp \in \Supp_R(M)\}.$$
Our main goal is to extend this theorem to modules finite over a
local homomorphism. To this end, we provide some preliminary
results. The most interesting one among them, is a generalization
of the Grothendieck's non-vanishing theorem (Theorem \ref{LC}).
\\
{\bf Theorem.} Let $(R,\fm)$ be a noetherian local ring and $M$ an
$R$-module finite over a local homomorphism. Then $n=dim_RM$ is
the greatest integer $i$ for which the $i$-th local cohomology
module $H^i_\fm(M)$ is non-zero.\\
Using the auxiliary results, we prove our main theorem (Theorem
\ref{mainthm}) with similar techniques as in \cite{KY} and
\cite{Y}. As a corollary, we also prove that if a local ring $R$
admits a module $M$ finite over a local homomorphism with finite
Gorenstein injective dimension and maximal dimension over $R$,
then $R$ is Cohen-Macaulay.

%
%

\noindent {\bf Convension.} Throughout this notes all rings are assumed to
be commutative,
notherian and local.
\section{Local cohomology over a local homomorphism}
Let $(R,\fm)$ be a noetherian local ring and $M$ a finite $R$-module of
the Krull dimension $n$. Then the Grothendieck's non-vanishing theorem
says that the local cohomology module $H^n_\fm(M)$ does not vanish (cf.
\cite[6.1.4]{BS}). In this section,
we prove the same non-vanishing result for a module finite over a
local homomorphism. First, we prove a couple of auxiliary lemmas.

The main tool we use in our proofs in this paper, is the
so-called {\it Cohen factorization} of a homomorphism.
\\
A local homomorphism $\varphi:R\to S$ is said to have a Cohen
factorization if it can be represented as a composition $R\to
R'\to S$ of local homomorphisms, where $R'$ is a complete local
ring, the map $R \to R'$ is flat with regular closed fiber $R'/\fm
R'$ and $R' \to S$ is a surjective homomorphism. For any local
homomorphism $\varphi:(R,\fm)\to (S,\fn)$ its semi-completion
$\varphi':R\to S\to \hat{S}$ admits a Cohen factorization (cf.
\cite{AFH}).

\begin{lem}\label{Supp}
Let $\varphi: (R,\fm) \to (S,\fn)$ be a flat local homomorphism.
If $M$ is an $S$-module, then
$$\begin{array}{c}
\Supp_R(M)=\{\varph^{-1}(\q) \, | \, \q \in \Supp_S(M) \} \, \,
\mathrm{and} \\
Supp_S(M)= \{ \fq \in \Spec(S)\, | \, \varphi^{-1}(\fq) \in
\Supp_R(M)\}.
\end{array}$$
\end{lem}
\begin{proof}
Suppose that $\fq \in \Spec(S)$ and set $\fp=\varph^{-1}(\q)$.
Therefore $\varphi$ induces a (faithfully) flat local homomorphism
$R_{\fp} \to S_{\fq}$ and the following isomorphisms hold.

$$\begin{array}{rcl}
  M_{\fp} \otimes_{R_{\fp}} S_{\fq} &\cong &(M \otimes_R R_{\fp})
\otimes_{R_{\fp}} S_{\fq} \\
                                    &\cong & M \otimes_R S_{\fq} \\
                                    &\cong & M \otimes_R (S \otimes_S
S_{\fq}) \\
                                    &\cong & M_{\fq} \otimes_R S
\end{array}$$
If $\fq \in \Supp_S(M)$ then, $\varphi$ being faithfully flat,
$M_{\fq} \otime_R S$ and  thus $M_{\fp} \otimes_{R_{\fp}} S_{\fq}$
are non-zero. Hence $\fp=\varph^{-1}(\q)\in \Supp_R(M)$.

On the other hand, if $\fp \in \Supp_R(M)$, using faithfully
flatness of $\varphi$ we can find a prime $\fq \in \Spec(S)$ such
that $\varph^{-1}(\q)= \fp$ and the above isomorphisms give $M_{\fq}
\otimes_R S \neq 0$ and hence $\fq \in \Supp_S(M)$.
\end{proof}

\begin{lem}\label{dim}
Let $\varphi: (R,\m) \to (S,\n)$ be a flat local homomorphism and
$M$ a non-zero $S$-module. Then $$\dim_SM \ge \dim_R M + \dim_S
S/{\fm S}.$$
\end{lem}
\begin{proof}
First note that since $\varphi$ is flat, every prime ideal of $S$
which contains $\fm S$ minimally, contracts to $\fm$ via
$\varphi$. Let $\fq \in \Spec(S)$ be one such prime. By Lemma
\ref{Supp}, $\fq \in \Supp_S(M)$.
\\
Suppose that $\dim_R M=n$ and $$\fp_0 \subset \cdots \subset
\fp_n=\fm$$ is a maximal chain of prime ideals in $\Supp_R(M)$.
Since $\varphi$ is flat and $\varphi^{-1}(\fq)=\fm$ we get a
strict chain in $\Spec(S)$, $$\fq_0 \subset \cdots \subset
\fq_n=\fq$$ such that $\varphi^{-1}(\fq_i)=\fp_i$ for $0 \leq i
\leq n$. Using Lemma \ref{Supp} again, this is a chain in
$\Supp_S(M)$ which gives the inequalities
$$\dim_S M \ge \dim_{S_{\fq}} M_{\fq} + \dim_S S/{\fq} \ge n +
\dim_S S/{\fq}.$$ In particular, we can choose $\fq$ such that
$\dim_S S/{\fm S}= \dim_S S/{\fq}$ to get the desired inequality.
\end{proof}

Now we are ready to prove the main theorem of this section.
\begin{thm}\label{LC}
Let $\varphi: (R,\m) \to (S,\n)$ be a local homomorphism and $M$ a
non-zero finite $S$-module. If $n= \dim_R M$ then
$$\H_{\fm}^n(M)\neq 0.$$
\end{thm}
\begin{proof}
First suppose that $\varphi$ is a flat homomorphism.
\\
As in the proof of Lemma \ref{dim}, we pick a prime $\fq \in
\Supp_S(M)$ which contains $\fm S$ minimally. Therefore, the
induced homomorphism $R \to S_{\fq}$ is faithfully flat and has
zero-dimensional fiber $S_{\fq}/{\fm S_{\fq}}$. Furthermore, the
inequality $\dim_{S_{\fq}} M_{\fq} \ge \dim_R M$ holds.

On the other hand, the isomorphism $\H_{\fm}^i(M_{\fq}) \cong
\H_{\fq S_{\fq}}^i(M_{\fq})$ holds for every $i \ge 0$. Since
$M_{\fq}$ is a finite $S_{\fq}$-module, $\H_{\fq
S_{\fq}}^d(M_{\fq}) \neq 0$ for $d = \dim_{S_{\fq}}M_{\fq}$ and
therefore we get $\dim_{S_{\fq}} M_{\fq} \leq \dim_R M_{\fq}$.

Note that $\Supp_R(M_{\fq}) \subseteq \Supp_R(M)$ and hence
$$\dim_{S_{\fq}} M_{\fq} \leq \dim_R M_{\fq} \leq \dim_R M \leq
\dim_{S_{\fq}} M_{\fq}.$$ Thus the equality $\dim_{S_{\fq}}
M_{\fq}=\dim_R M=n$ holds.

Finally, the following isomorphisms give the desired non-vanishing
in this case. $$\H_{\fm}^n(M_{\fq}) \cong \H_{\fm
S_{\fq}}^n(M_{\fq}) \cong \H_{\fm}^n(M) \otimes_R S_{\fq}$$
In general, $R \to S \to \hat{S}$ admits a Cohen factorization $R
\to R' \to \hat{S}$. Since $M$ is finite over $S$, so is
$\hat{M}=M \otimes_S \hat{S}$ as an $\hat{S}$-module and thus as
an $R'$-module. It is also easy to see that
$\Supp_R(\hat{M})=\Supp_R(M)$ and consequently we have $\dim_R
\hat{M}= \dim_R M =n$. So using the first part of the proof, we
get
$$\H_{\fm}^n(\hat{M}) \neq 0.$$ Using the isomorphisms $$\H_{\fm}^n(\hat{M}) \cong \H_{\fm
\hat{S}}^n(\hat{M}) \cong \H_{\fm S}^n(M) \otimes_S \hat{S}$$ we
conclude that $\H_{\fm S}^n(M)\cong \H_{\fm}^n(M)$ is non-zero.

\end{proof}
\section{A Generalized Bass Formula}
This section is dedicated to proving a Bass-type formula for a module
which
is finite over a local homomorphism and has finite Gorenstein injective
dimension.
\begin{defn}
An $R$-module $G$ is said to be {\rm Gorenstein injective} if and
only if there exists an exact complex of injective $R$-modules,
$$I=\cdots \to I_2\lo I_1\lo I_0\lo I_{-1}\lo
I_{-2}\lo\cdots$$ such that the complex $\Hom_R(J,I)$ is exact for
every injective $R$-module $J$ and $G$ is the kernel in degree 0
of $I$. The {\rm Gorenstein injective dimension} of an $R$-module
$M$, $\Gid_R (M)$, is defined to be the infemum of integers $n$
such that there exists an exact sequence $$0 \to M \to G_0 \to
G_{-1} \to \cdots \to G_{-n} \to 0$$ with all $G_i$'s Gorenstein
injective.
\end{defn}

\begin{prop}\label{width}
Let $\varphi: (R, \fm) \to (S, \fn)$ be a local homomorphism and
$M$ a finite $S$-module. If $\fp \in \Supp_R(M)$ then
$M_{\fp}/{\fp M_{\fp}} \neq 0$.
\end{prop}
\begin{proof}
First suppose that $\varphi$ admits a Cohen factorization $R\to
R'\to S$. Since $R' \to S$ is surjective, $M$ is also finite as an
$R'$-module. Thus we can assume that $\varphi$ is flat. By Lemma
\ref{Supp}, there exists a $\fq \in \Supp_S(M)$ such that
$\varphi^{-1}(\fq)=\fp$.
\\
By NAK we have $M_{\fq}/{\fp M_{\fq}} \neq 0$. Since $M_{\fq}/{\fp
M_{\fq}} \cong(T^{-1}M/\fp T^{-1}M)_{T^{-1}\fq},$ where
$T=\varphi(R-\fp),$ we get $M_{\fp}/{\fp M_{\fp}} \cong
T^{-1}M/\fp T^{-1}M \neq 0$.
\\
In general $R \to S \to \hat{S}$ admits a Cohen factorization. If
$\fp \in \Supp_R(M)$ then $\fp \in \Supp_R(M \otimes_S \hat{S})$
as well and then, setting $\hat{M}=M \otimes_S \hat{S}$, we have
$M_{\fp}/{\fp M_{\fp}} \otimes_S \hat{S}
\cong \hat{M}_ {\fp}/{\fp \hat{M}_{\fp}} \neq 0$ and we are done.
\end{proof}

\begin{thm}\label{mainthm}
Let $R$ be a noetherian local ring and $M$ be an $R$-module finite
over a local homomorphism with $\gid_R(M)<\infty$. Then $$\dim_R M
\leq \gid_R(M)=\sup\{\depth R_\fp \, | \, \fp\in\Supp_R(M)\}.$$
\end{thm}
\begin{proof}
To prove the first inequality we use the proof of \cite[3.4]{S} to
see that $\H_\fm^i(M)=0$ for $i> \Gid_R(M)$. Thus, by Theorem
\ref{LC}, we have $$\dim_R M \leq \Gid_R(M).$$



To prove the formula, we use \cite[2.18]{CFH}, to get an exact
sequence $$0\to K\to L\to M\to 0$$ with $K$ Gorenstein injective
and $\id_R(L)=\gid_R(M)$. By \cite{Chouinard}, we have
$$\id_R(L)=\sup\{\depth R_\fp-\width_{R_\fp}
L_\fp \, | \,
\fp\in\Supp_R(L)\},$$ where $\width_{R_\fp}M_\fp = \inf \{ \,i \,
| \Tor_i^{R_\fp}(R_\fp/\fp R_\fp,M_\fp)\neq0 \}$.
\par
\bigskip
For any $\fp\in\Supp_R(M)$, we
get an exact sequence $L_\fp/{\fp L_\fp} \to M_\fp/{\fp M_\fp} \to
0$ which gives $L_\fp/\fp L_\fp\neq 0$, by Proposition
\ref{width}. Therefore we have
$$\begin{array}{c}
\sup \{ \depth R_\fp - \width_{R_\fp} L_\fp \, | \, \fp \in
\Supp_R(L) \} \geq \\
\sup \{ \depth R_\fp - \width_{R_\fp} L_\fp \, | \, \fp \in
\Supp_R(M) \}= \\
\sup \{\depth R_\fp \, | \, \fp \in \Supp_R(M)\} \geq 0.
\end{array}$$
This proves the desired formula for $M$ a Gorenstein injective
$R$-module.
\par
For any $\fp\in\Supp_R(L)-\Supp_R(M)$, we get
$L_\fp\cong K_\fp$. Since $K$ is Gorenstein injective, there
exists an exact sequence
$$\cdots\to I_1\to I_0\to L_\fp\to 0,$$ where for each $\ell\ge
0$, $I_\ell$ is an injective $R_\fp$-module. \\ Set $K_\ell= \ker(
I_{\ell-1}\to I_\ell)$. Then for any $R_{\fp}$-module $T$ with
$\pd_{R_\fp}(T)=t < \infty$ we have $\Ext^i_{R_\fp}(T,L_\fp)\cong
Ext^{i+t}_{R_\fp}(T,K_t)=0$ for $i>0$. Thus, by
\cite[5.3(c)]{CFF}, we get
$$\depth R_\fp -\width_{R_\fp}L_\fp \leq 0.$$ \\
Therefore, if $\gid_R(M)=\id_R(L)>0$ then we have

$$\begin{array}{c}
\Gid_R(M)=\id_R(L)=
\\
\sup \{ \depth R_\fp -\width_{R_\fp}L_\fp \, | \, \fp \in
\Supp_R(L)\}= \\
\sup \{ \depth R_\fp -\width_{R_\fp}L_\fp \, | \, \fp \in
\Supp_R(M)\}= \\
\sup \{ \depth R_\fp \, | \, \fp \in \Supp_R(M)\}.
\end{array}$$
\end{proof}
The following corollary of Theorem \ref{mainthm}, gives a
sufficient condition for Cohen-Macaulayness of a noetherian local
ring. This result generalizes \cite [3.5]{T} as well as
\cite[1.3]{Y}.
\begin{cor}\label{cmness}
Let $R$ be a local ring and $M$ an $R$-module finite over a local
homomorphism. If $M$ has finite Gorenstein dimension and maximal
Krull dimension over $R$ (i.e. $\Gid_R(M) < \infty$ and $\dim_R M
=\dim R$), then $R$ is Cohen-Macaulay.
\end{cor}
\begin{proof}
By \ref{mainthm}, there is a prime ideal $\fp \in \Supp_R(M)$ such
that $$ \dim R = \dim_R M \leq \Gid_R(M)= \depth R_\fp \leq \dim
R_\fp \leq \dim R.$$ Therefore, $\fp$ must be the maximal ideal of
$R$ and thus $R$ is Cohen-Macaulay .
\end{proof}


\begin{cor}\label{cb}
Let $R$ be an almost Cohen-Macaulay local ring (i.e.
$\dim(R)-\depth(R)\le 1$) and let $M$ be an $R$-module finite over
a local homomorphism. If $M$ has finite Gorenstein injective dimension
over $R$ then $\gid_R(M)=\depth R$.
\end{cor}

\begin{proof}
It is enough to use Theorem \ref{mainthm} and the fact that over
an almost Cohen-Macaulay ring $R$ the inequality $\depth R_\fp
\leq \depth R_\fq$ holds for $\fp \subseteq \fq$ in $\Spec(R)$.
\end{proof}
Note that Corollary \ref{cb} can be considered as a generalized
Bass formula for Gorenstein injective dimension.

%
%
We conclude this paper with a change of ring result for Gorenstein
injective dimension .

\begin{cor}
Let $(R, \fm)$ be a local ring and $M$ be an $R$-module finite over a
local homomorphism. If
$x\in\fm$ is an $R$- and $M$-regular element, then
$$\gid_{R/xR}(M/xM)\le\gid_R(M)-1.$$
Furthermore, the equality holds when $R$ is almost Cohen-Macaulay
and $\gid_R(M)< \infty.$
\end{cor}
\begin{proof}
If $M$ has finite Gorenstein injective dimension over $R$ then the
proof of \cite[Lemma 2]{SSY} shows that $\gid_{R/xR}(M/xM)< \infty,$
too . By \ref{width}, for any $\fp \in \Supp_R(M)$ with $x \in
\fp$, the $R_\fp$-module $M_\fp/xM_\fp$ is non-zero. Therefore,
$$\Supp_{R/xR}(M/xM)=\{\fp/xR | \fp \in \Supp_R(M) \, \,
\mathrm{with} \, \, x \in \fp \}.$$ Thus Theorem \ref{mainthm}
gives the desired inequality. In the case of an almost
Cohen-Macaulay base ring the equality is a consequence of
\ref{cb}.
\end{proof}
\providecommand{\bysame}{\leavevmode\hbox
to3em{\hrulefill}\thinspace}

\end{document}